\documentclass[12pt]{amsart}
\newtheorem{theorem}{Th\'eor\`eme}

 \newtheorem{corollary}{Corollaire}

\usepackage[latin1]{inputenc}

\input{amssym.def}

\begin{document}

\title{Distances invariantes et points fixes d'applications holomorphes} 

\author{Jean-Pierre Vigu\'e}
\maketitle
\begin{abstract}
In this paper, we prove the following result : let $X$ be a complex manifold, hyperbolic for
the Carathéodory distance and let $U$ be an open set relatively compact in $X$. Then, 
there exists $k<1$ such that we get, for the Carathéodory infinitesimal metric 
$E_X(x,v)\leq kE_U(x,v)$. We also get results concerning fixed points of holomorphic 
mappings from $X$ to $U$.
\bigskip

\noindent{\sc Résumé.}
Dans cet article, nous montrons le résultat suivant : Soit $X$ une variété complexe hyperbolique
pour la distance de Carathéodory et soit $U$ un ouvert relativement compact dans $X$. Alors, 
il  existe $k<1$ tel que, pour la métrique infinitésimale de Carathéodory, on ait  
$E_X(x,v)\leq kE_U(x,v)$. Nous obtenons aussi des résultats sur les points fixes d'applications 
holomorphes de $X$ dans $U$.

\bigskip

\noindent{\bf Mathematics subject classification} (2010)
32F45, 32H02.
\bigskip

\end{abstract}

\section{Introduction}

Un certain nombre d'auteurs ont considéré le problème suivant : soit 
$X$ une variété analytique complexe munie d'une distance invariante $d_X$. 
Soit $U$ un ouvert relativement compact dans $X$. Existe-t-il une constante 
$k<1$ telle que, pour tout $x\in U$, pour tout $y\in U$, on ait
$$d_{X}(x,y)\leq kd_{U}(x,y).$$
Si c'est le cas, on en déduit que, pour toute application holomorphe 
 $ f:X \longrightarrow X$ telle que $f(X)$ 
soit relativement compact dans $X$, $f$ admet un point fixe unique 
$a\in X$, et $a=\lim f^{n}(x_{0})$, 
où $x_0$ est un point quelconque de $X$.

Sur cette question, on peut citer en particulier les deux résultats suivants :

- H. Reiffen \cite{re} annonce ce résultat dans le cas d'un espace analytique complexe 
$X$ $c_X^i$-hyperbolique (ce qui signifie que la pseudodistance intégrée de Carathéodory $c_X^i$
 est une distance sur $X$).

- C. Earle et R. Hamilton \cite{ea} montrent ce résultat dans le cas d'un domaine 
borné $X$ d'un espace de Banach complexe et d'un ouvert $U\subset X$ complètement intérieur à $X$, 
ce qui signifie que la distance de $U$ à la frontière de $X$ est strictement positive. 
Dans le cas de la dimension finie, ceci signifie seulement que $U$ est relativement compact dans $X$.

Il nous faut aussi signaler un résultat de M. Hervé \cite{he}, p. 92 sur les points
 fixes d'applications holomorphes, résultat obtenu sans utiliser les distances invariantes : soit 
$U$ un ouvert de $\Bbb C^n$ et soit $ f:U \longrightarrow K$ une application holomorphe,
 où $K$ est un compact de $U$. Alors 
$f$ admet un point fixe unique.

Enfin, signalons un résultat de P. Mazet et J.-P. Vigué \cite{ma} (voir aussi \cite{vi}) obtenu en utilisant les résultats 
précédents : l'ensemble des points fixes d'une application holomorphe $ f:X \longrightarrow X$, où $X$ 
est un domaine borné de $\Bbb C^n$, est une sous-variété analytique complexe de $X$.

Nous allons maintenant montrer deux résultats sur cette question, un obtenu en utilisant la distance 
intégrée de Carathéodory, l'autre en utilisant la distance de Kobayashi. Pour commencer, nous allons 
rappeler quelques propriétés des distances invariantes.

[Je remercie Jean-Jacques Loeb pour toutes les discussions que nous avons eues
 sur cette question.]

\section{Rappel sur les distances et métriques invariantes} 
Soit $X$ une variété analytique complexe. 
On d\'efinit la pseudom\'etrique infinit\'esimale
 de Carath\'eodory (ou de Reiffen-Cara\-th\'eo\-dory) de la fa\c con suivante
 : $\forall x\in X, \forall v\in T_{x}X$,  
$$E_{X}(x,v)= \sup |\varphi' (x).v|,$$
où $\varphi$ parcourt l'ensemble des applications  holomorphes de $X$ 
dans le disque-unité $\Delta $.
On v\'erifie facilement que $ E_{X}
 $ est une pseudom\'etrique invariante dans le sens suivant : si $ f:X\longrightarrow X' $ 
est une application holomorphe, on a,
$\forall x\in X, \forall v\in T_{x}X$,   
$$E_{X'}(f(x),f'(x).v)\leq E_{X}(x,v).$$
 En particulier, les isomorphismes analytiques sont des isom\'etries.

De mani\`ere duale, on d\'efinit la pseudom\'etrique infinit\'esimale
 de Kobayashi (ou de Royden-Kobayashi):
$\forall x\in X, \forall v\in T_{x}X$,
$$F_{X}(x,v)=\inf_{\varphi\in H(\Delta,X),\varphi(0)=x}\{|\lambda|
 \hbox{\rm\ tel que }\lambda\varphi' (0)=v\}.$$
 On v\'erifie que $ F_{X} $ est elle aussi une pseudom\'etrique invariante.

Etant donnée une pseudométrique invariante sur une variété $X$, on peut définir par intégration la 
longueur d'un chemin de classe $C^1$ par morceaux. Ensuite, on définit une pseudodistance intégrée. La 
pseudodistance de deux points $a$ et $b$ est définie comme la borne inférieure de la longueur des chemins 
d'origine $a$ et d'extrémité $b$ (voir \cite{ko}). C'est, en un sens facile à préciser, une pseudodistance 
invariante.
Par cette technique, en partant de la pseudométrique infinitésimale de Carathéodory, on définit 
la pseudodistance intégrée de Carathéodory $c_X^i$ ;  en partant de la pseudométrique infinitésimale
 de Kobayashi, on définit 
la pseudodistance de Kobayashi $k_X$. Par exemple, dans le cas du disque-unité $\Delta$ dans $\Bbb C$, on obtient
$$c_\Delta^i(z,w)=k_\Delta(z,w)=\omega (z,w)= \tanh ^{-1}\Big|\,{z-w      \over 1-\overline{w}z
     }\,\Big|,$$
 où $\omega $ est la distance de Poincaré.

Enfin, on d\'efinit la pseudodistance de Carath\'eodory $ c_{X} $ sur une
variété $X$ de la fa\c con suivante : pour tous $ a $ et $ b $ dans $ X$, 
$$c_{X}(a,b)=\sup_{\varphi\in H(X,\Delta)} \omega (\varphi(a),\varphi(b)),$$
 o\`u $ H(X,\Delta) $ d\'esigne l'ensemble des applications 
holomorphes de $ X $ dans
 le disque-unit\'e $ \Delta$. Bien sûr, $c_{X}$ est une pseudodistance invariante. 

Pour une description plus détaillée des pseudométriques et pseudodistances invariantes, voir 
le livre de S. Kobayashi \cite{ko}.

\section{Utilisation de la distance intégrée de Carathéodory}

Soit $X$ une variété analytique complexe et soit $c_X^i$ la pseudodistance de 
Carathéodory sur $X$. On dit que $X$ est 	$c_X^i$-hyperbolique si $c_X^i$ est une distance sur $X$.
Si $c_X^i$ est une distance, on déduit du fait que c'est une distance intégrée et 
d'un résultat classique qu'elle définit la 
topologie de $X$. Nous avons le théorème suivant.

\begin{theorem}
Soit $X$ une variété analytique complexe $c_X^i$-hyperbolique et soit $U\subset X$ un domaine 
relativement compact dans $X$. Alors, il existe une constante $k<1$ telle que,
 $\forall a\in X, \forall v\in T_{a}X$, on ait 

$$E_{X}(a,v)\leq kE_{U}(a,v).$$ 
Plus précisément, $k=\tanh M$, où $M$ est le diamètre de $U$ pour
 la pseudodistance de Carathéodory.
 
Par intégration, on en déduit que, pour tous $a$ et $b$ appartenant à $U$, 
$$c_X^i(a,b)\leq kc_U^i(a,b)$$
\end{theorem}

{\it Démonstration.}  Soit $c_X$ la pseudodistance de Carathéodory sur $X$. Comme 
$\overline{U}$  est compact, le diamètre de $U$ pour la pseudodistance de Carathéodory
 qui, par définition, vaut 
$ \sup_{x\in\overline{U},y\in\overline{U}}c_{X}(x,y)$
est une constante $M<+\infty $, et
$$M=\sup_{f\in H(X,\Delta),f(x)=0}\omega(0,f(y))=\sup_{f\in
 H(X,\Delta),f(x)=0}\tanh^{-1}\Big|f(y)\Big|.$$On en déduit immédiatement que, si on pose $k=\tanh M$, on a : 
$ k<1 $ et, $\forall a\in X, \forall f\in H(X,\Delta)$, tel que $ f(a)=0$,
 on a : $ \|f\|_{\overline{U}} \leq k<1$.

Etant donn\'ee une fonction holomorphe $ f:X\longrightarrow\Delta
 $ telle que $ f(a)=0$, on d\'efinit une fonction $ \varphi:U\longrightarrow\Delta $ 
par la formule
$$\varphi(z)=(1/k)f(z).$$ 
Il est clair que $ \varphi $ envoie $ U $ dans $ \Delta$, et on a : 

$$\Big|\varphi' (z).v\Big|=(1/k)\Big|f' (z).v\Big|.$$
En prenant la borne sup\'erieure, pour toutes les fonctions holomorphes
 $ f:X\longrightarrow\Delta$, on en d\'eduit que 
$$E_{U}(a,v)\geq(1/k)E_{X}(a,v),$$
 ce qui donne 
$$E_{X}(a,v)\leq kE_{U}(a,v),$$
ce qui est le r\'esultat annonc\'e. Par int\'egration, on en d\'eduit
 que, pour tout $ a\in U$, pour tout $ b\in U$, $ c^{i}_{X}(a,b)\leq
 kc^{i}_{U}(a,b)$. 

Comme application, nous avons le th\'eor\`eme suivant. 
\begin{theorem}
Soit $ X $ une vari\'et\'e $c_X^i$-hyperbolique et soit $ f:X\longrightarrow
 X $ une application holomorphe telle que $ f(X) $ soit relativement compact
 dans $ X$. Alors la suite des it\'er\'ees $ f^{n} $ converge, pour la
 topologie compacte ouverte vers une application constante $ z\mapsto c$,
 o\`u $ c $ est l'unique point fixe de $ f$. 

(En particulier, l'hypoth\`ese que $ X $ est $c_X^i$-hyperbolique est v\'erifi\'ee
 si on suppose que la pseudodistance de Carath\'eodory $ c_{X}$ sur $ X
 $ est une distance sur $ X$.)
\end{theorem}

{\it Démonstration.}
On d\'eduit facilement du fait que $ f(X) $ est relativement compact
 dans $ X$ qu'il existe un ouvert $ U $ relativement compact dans $ X $ tel
 que $ f(X)\subset U\subset X$. Comme $ f $ envoie $ X $ dans $ U$, on
 a, pour tout $ a\in X$, pour tout $ v\in T_{a}(X)$, 
$$E_{U}(f(a),f'(a).v)\leq E_{X}(a,v).$$
 D'apr\`es le th\'eor\`eme 1, il existe une constante $ k<1 $ telle
 que 
$$E_{X}(f(a),f'(a).v)\leq kE_{U}(f(a),f'(a).v).$$
 On trouve alors 
$$E_{X}(f(a),f'(a).v)\leq kE_{X}(a,v).$$
 Par int\'egration, on en d\'eduit que, pour tout $ a\in X$, pour tout
 $ b\in X$, 
$$c_X^i(f(a),f(b))\leq kc_X^i(a,b).$$
 En utilisant ce r\'esultat pour $ b=f(a)$, et en it\'erant, on en
 d\'eduit que, pour tout $ n>0$, 
$$c_X^i(f^{n}(a),f^{n+1}(a))\leq k^{n}c_X^i(a,f(a)).$$
 De fa\c con tout \`a fait classique, on en d\'eduit que, pour tout
 $ a\in X$, la suite des it\'er\'ees $ (f^{n}(a)) $ est une suite de Cauchy
 pour $c_X^i$ sur $ \overline{f(X)} $ qui est compact. Elle converge donc
 vers $ c $ qui est l'unique point fixe de $ f$. 

Maintenant, soit $ K $ un compact de $ X$. On a : 
$$\sup_{z\in K}c_X^i(z,f(z))=M<+\infty.$$ 
Pour $ n\geq p$, on a pour tout $ z\in K$, 
$$c_X^i(f^{n}(z),f^{p}(z))\leq k^{p-1}/(1-k)M,$$ 
ce qui montre la convergence uniforme sur $ K $ de la suite $ (f^{n}(z))
 $ vers $ c$. Ainsi, la suite $ (f^{n}(z)) $ converge vers la fonction
 constante \'egale \`a $ c $ pour la topologie compacte ouverte. 

Nous avons aussi le corollaire suivant.
 
\begin{corollary}
Soit $ A $ une vari\'et\'e de Stein, et soit $ X\subset A $ un domaine 
relativement compact dans $A$.
 Soit $ U $ un ouvert de $ X $ relativement compact dans $ X $ et soit 
$ f:X\longrightarrow U $ une application holomorphe. Alors $ f $ admet
 un point fixe unique $ c\in U $ et la suite des it\'er\'ees $ (f^{n})
 $ converge vers $ c $ pour la topologie compacte ouverte.  
\end{corollary}

{\it Démonstration.} D'après l'hypothèse, les fonctions holomorphes sur $A$ 
séparent les points de $X$. Comme $\overline {X}$ est compact, ces fonctions 
sont bornées sur $X$. Par suite, $X$ est $c_X$-hyperbolique et $c_X^i$-hyperbolique, 
et on peut appliquer le théorème 2. 

Remarquons que l'on peut aussi appliquer ce r\'esultat pour des vari\'et\'es
 qui ne sont pas hyperboliques pour la pseudodistance int\'egr\'ee
 de Ca\-rath\'eodory. Plus pr\'ecis\'ement, nous avons le corollaire suivant.
 
\begin{corollary}
Soit $ A $ une vari\'et\'e de Stein, et soit $ X\subset A $ un domaine.
 Soit $ U $ un ouvert de $ X $ relativement compact dans $ X $ et soit 
$ f:X\longrightarrow U $ une application holomorphe. Alors $ f $ admet
 un point fixe unique $ c\in U $ et la suite des it\'er\'ees $ (f^{n})
 $ converge vers $ c $ pour la topologie compacte ouverte. 
\end{corollary}
Par un raisonnement classique de compacit\'e, on construit un ouvert
 $ U' $ relativement compact dans $ X $ tel que $ U $ soit relativement 
compact dans $U'$. On remarque
 alors que, comme $ A $ est une vari\'et\'e de Stein, les fonctions holomorphes
 sur $ A $ s\'eparent les points de $ U' $, et comme $ U'  $ est relativement
 compact dans $ X$, les restriction \`a $ U' $ des fonctions holomorphes
 sur $ A $ sont born\'ees sur $ U' $, ce qui entra\^{\i}ne que $ U' $ est
 hyperbolique pour la distance de Carath\'eodory. D'apr\`es le th\'eor\`eme
 2, $ f $ admet un point fixe unique $ c\in U$. 

Le fait que la suite $ (f^{n}) $ converge vers $ c $ pour la topologie
 compacte ouverte se montre de la fa\c con suivante : soit $ K $ un compact
 contenu dans $ X$. Quitte \`a changer $ U' $, on peut supposer que $ U'
 $ contient $\overline {U}$ et $ K$. Il suffit alors d'appliquer le théorème 2. 

Pour terminer ce paragraphe, remarquons que le th\'eor\`eme que nous
 venons de d\'emontrer permet de retrouver le r\'esultat de M. Herv\'e \cite{he}.
 En effet, $\Bbb C^n$ est une vari\'et\'e de Stein et on peut appliquer le
 corollaire pr\'ec\'edent.  

\section {Utilisation de la distance de Kobayashi }

On peut se demander si on peut g\'en\'eraliser le th\'eor\`eme 1 en
 rempla\c cant la m\'etrique infinit\'esimale de Carath\'eodory par
 la m\'etrique infinit\'esimale de Kobayashi. Rappelons qu'une vari\'et\'e
 complexe $ X $ est hyperbolique si la pseudodistance de Kobayashi 
$ k_{X} $ est une distance sur $ X$. La premi\`ere remarque est qu'il existe
 des vari\'et\'es hyperboliques compactes. Dans ce cas, si on consid\`ere
 $ U=X\subset X$, on ne peut pas esp\'erer obtenir une in\'egalit\'e
 stricte entre $ F_{U}(a,v) $ et $ F_{X}(a,v)$. Pour obtenir un r\'esultat
 de ce genre, il faut faire des hypoth\`eses suppl\'ementaires sur
 les vari\'et\'es consid\'er\'ees (afin d'\'eliminer le cas compact).
 Plus pr\'ecis\'ement, nous avons le th\'eor\`eme suivant.  
\begin{theorem}
Soit $ A $ une vari\'et\'e de Stein de dimension $ n$. Soit $ X\subset
 A $ un domaine hyperbolique et soit $ U $ un ouvert de $ X $ relativement
 compact dans $ X$. Alors, il existe une constante $ k<1 $ telle que,
 pour tout $ a\in U$, pour tout $ v\in T_{a}U$, on ait 
$$F_{X}(a,v)\leq kF_{U}(a,v).$$
 De m\^eme, pour tout $ a\in U$, pour tout $ b\in U$, on a 
$$k_{X}(a,b)\leq kk_{U}(a,b).$$
 \end{theorem}
{\it Démonstration.}
D'apr\`es la th\'eorie des vari\'et\'es de Stein, il existe 
 un plongement ferm\'e $ i $ de $ A $ sur une vari\'et\'e ferm\'ee $ i(A)
 $ de $\Bbb C^{2n+1}$.

On d\'eduit alors du fait que $ i(U) $ est relativement compact dans
$\Bbb C^{2n+1}$ qu'il existe $ R>0 $ tel que, $\forall x\in i(U)$, $ i(U)\subset
 B(x,R)$, o\`u $ B(x,R) $ est la boule de centre $ x $ et de rayon R pour
 une norme quelconque sur $\Bbb C^{2n+1}$. 

D'apr\`es R. Gunning et H. Rossi \cite{gu}, chapitre VIII, section C,
 th\'eor\`eme 8, il existe un voisinage $ W $ de $ i(A) $ et une r\'etraction
 holomorphe $ \rho $ de $ W $ sur $ i(A)$. D'autre part, quitte \`a diminuer
 la taille de $ W$, on peut supposer qu'il existe  $ r>0 $ tel que 
$$\forall x\in i(U), B(x,r)\subset W,$$
$$\forall x\in i(U), \rho(B(x,r))\subset i(X).$$

Soit maintenant $ \varphi:\Delta\longrightarrow i(U) $ une
 application holomorphe telle que $ \varphi(0)=a$. On peut alors d\'efinir
 une application holomorphe $ \psi:\Delta\longrightarrow
 i(X) $ par la formule 
$$\psi(\zeta)=\rho((1+r/R)(\varphi(\zeta)-\varphi(0))+\varphi(0)).$$
 Il faut d'abord v\'erifier que $ \psi $ est bien d\'efinie. Pour cela,
 on \'ecrit 
$$(1+r/R)(\varphi(\zeta)-\varphi(0))+\varphi(0)=\varphi(\zeta)+r/R(\varphi(\zeta)-\varphi(0)).$$
 Or, $ \varphi(\zeta)\in i(U), \|\varphi(\zeta)-\varphi(0)\|<R, \|r/R(\varphi(\zeta)-\varphi(0))\|<r$,
 et $ W $ contient la boule de centre $ \varphi(\zeta) $ et de rayon $ r$.
 Par suite, $ (1+r/R)(\varphi(\zeta)-\varphi(0))+\varphi(0)\in W$, et
 son image par $ \rho $ est contenue dans $ i(X)$. 

D'autre part, il est clair que $ \psi(0)=a $ et $ \psi'(0)=(1+r/R)\varphi'(0)$.
 Ceci suffit \`a d\'emontrer que 
$$F_{X}(a,v)\leq1/(1+r/R)F_{U}(a,v).$$
 Par int\'egration, on en d\'eduit la m\^eme in\'egalit\'e sur la
 distance de Kobayashi. 

Par les m\^emes m\'ethodes que dans la section 3, on montre un th\'eor\`eme
 sur les points fixes d'applications holomorphes. 
\begin{theorem}
Soit $ A $ une vari\'et\'e de Stein. Soit $ X\subset
 A $ un domaine hyperbolique et soit $ f:X\longrightarrow
 X $ une application holomorphe telle que $ f(X) $ soit relativement compact
 dans $ X$. Alors la suite des it\'er\'ees $ f^{n} $ converge, pour la
 topologie compacte ouverte vers une application constante $ z\mapsto c$,
 o\`u $ c $ est l'unique point fixe de $ f$.  
\end{theorem}
A partir de ces consid\'erations, on peut donner une nouvelle d\'e\-mons\-tra\-tion
 du corollaire 2 qui utilise la m\'etrique infinit\'esimale de Kobayashi.
 Le r\'esultat obtenu est, il semble, tout \`a fait semblable \`a celui
 du corollaire 2.

\bigskip\bigskip

 LMA, Universit\'{e} de Poitiers, CNRS, 
Math\'{e}matiques, SP2MI, BP 30179, 86962 FUTUROSCOPE.
\smallskip

e-mail :  vigue@math.univ-poitiers.fr ou jp.vigue@orange.fr

\end{document}